\documentclass[12pt]{article}

\usepackage[english]{babel}
\usepackage[utf8]{inputenc}
\usepackage[T1]{fontenc}
\usepackage{xcolor}
\usepackage{mathtools}
\usepackage{lipsum}
\usepackage{multicol}
\makeatletter
\newcounter {subsubsubsection}[subsubsection]
\renewcommand\thesubsubsubsection{\thesubsubsection .\@arabic\c@subsubsubsection}
\newcommand\subsubsubsection{\@startsection{subsubsubsection}{4}{\z@}%
                                     {-3.25ex\@plus -1ex \@minus -.2ex}%
                                     {1.5ex \@plus .2ex}%
                                     {\normalfont\normalsize\bfseries}}
\renewcommand\paragraph{\@startsection{paragraph}{5}{\z@}%
                                    {3.25ex \@plus1ex \@minus.2ex}%
                                    {-1em}%
                                    {\normalfont\normalsize\bfseries}}
\renewcommand\subparagraph{\@startsection{subparagraph}{6}{\parindent}%
                                       {3.25ex \@plus1ex \@minus .2ex}%
                                       {-1em}%
                                      {\normalfont\normalsize\bfseries}}
\newcommand*\l@subsubsubsection{\@dottedtocline{4}{10.0em}{4.1em}}
\renewcommand*\l@paragraph{\@dottedtocline{5}{10em}{5em}}
\renewcommand*\l@subparagraph{\@dottedtocline{6}{12em}{6em}}
\newcommand*{\subsubsubsectionmark}[1]{}
\makeatother

\usepackage{hyperref}
\newcommand\restr[2]{{
  \left.\kern-\nulldelimiterspace 
  #1 
  \vphantom{\big|} 
  \right|_{#2} 
  }}
\makeatletter
\def\toclevel@subsubsubsection{4}
\def\toclevel@paragraph{5}
\def\toclevel@subparagraph{6}
\makeatother

\usepackage[top=1.3cm, bottom=2.0cm, outer=2.5cm, inner=2.5cm, heightrounded,
marginparwidth=1.5cm, marginparsep=0.4cm, margin=2.5cm]{geometry}

\usepackage{graphicx} 
\usepackage{amsmath}  
\usepackage{amsfonts} %
\usepackage{amssymb}  %
\usepackage{amsthm}
\usepackage{bm}
\usepackage{hyperref}
\hypersetup{
    colorlinks = false,
    linkbordercolor = {red},
}

\usepackage[
backend=biber,
style=alphabetic,
sorting=anyvt,
natbib=true
]{biblatex}
\usepackage{amsmath,amsthm,amssymb}
\usepackage{enumitem}

\newcommand{\C}{\mathbb{C}}

\DeclareFontFamily{OT1}{pzc}{}
\DeclareFontShape{OT1}{pzc}{m}{it}{<-> s * [1.10] pzcmi7t}{}
\DeclareMathAlphabet{\mathpzc}{OT1}{pzc}{m}{it}

\newcommand{\overbar}[1]{\mkern 1.5mu\overline{\mkern-1.5mu#1\mkern-1.5mu}\mkern 1.5mu}
\newcommand{\norm}[1]{\left\lVert#1\right\rVert}
\newcommand{\abs}[1]{\left\lvert#1\right\rvert}

\newtheorem{theorem}{Theorem}
[section]
\newtheorem{thmx}{Theorem}

\theoremstyle{definition}
\newtheorem{definition}{Definition}[section]

\theoremstyle{plain}

\newtheorem{proposition}[theorem]{Proposition}

\theoremstyle{remark}
\newtheorem{remark}[theorem]{Remark}

\theoremstyle{definition}
\newtheorem{example}[theorem]{Example}

\newcommand{\I}{\sqrt{-1}}
\newcommand{\del}{\partial}
\newcommand{\dbar}{\bar{\partial}}
\newcommand{\ddbar}{\del\dbar}

\newcommand{\SBG}{\mathcal{S}\mathcal{B}\mathcal{G}_1}
\newcommand{\SBGc}{\mathcal{S}\mathcal{B}\mathcal{G}_c}

\newcommand{\ricci}{\mathrm{\mathbf{Ric}}}

\newcommand{\holosections}{\Gamma_{\mathcal{O}}}

\newcommand{\HH}{\mathcal{H}^2}
\usepackage{tikz}
\usetikzlibrary{shapes.geometric, arrows}
\tikzstyle{startstop} = [rectangle, rounded corners, minimum width=3cm, minimum height=1cm,text centered, draw=black]
\tikzstyle{process} = [rectangle, minimum width=3cm, minimum height=1cm, text centered, draw=black]

\tikzstyle{arrow} = [thick,->,>=stealth]

\makeatletter
\long\def\pgfshapeaddanchor#1#2{%
{%
  \def\pgf@sm@shape@name{#1}%
  \let\anchor=\pgf@sh@anchor%
  #2}%
}
\pgfshapeaddanchor{rectangle}{%
  \anchor{west north west}{%
    \pgf@process{\northeast}%
    \pgf@ya=1.5\pgf@y%
    \pgf@process{\southwest}%
    \pgf@y.5\pgf@y%
    \advance\pgf@y by \pgf@ya%
    \pgf@y=.5\pgf@y%
  }%
  \anchor{north north west}{%
    \pgf@process{\southwest}%
    \pgf@xa=1.5\pgf@x%
    \pgf@process{\northeast}%
    \pgf@x.5\pgf@x%
    \advance\pgf@x by \pgf@xa%
    \pgf@x=.5\pgf@x%
  }%
  \anchor{west south  west}{%
    \pgf@process{\northeast}%
    \pgf@ya=.5\pgf@y%
    \pgf@process{\southwest}%
    \pgf@y=1.5\pgf@y%
    \advance\pgf@y by \pgf@ya%
    \pgf@y=.5\pgf@y%
  }%
  \anchor{south south  west}{%
    \pgf@process{\northeast}%
    \pgf@xa=.5\pgf@x%
    \pgf@process{\southwest}%
    \pgf@x=1.5\pgf@x%
    \advance\pgf@x by \pgf@xa%
    \pgf@x=.5\pgf@x%
  }%
}

\makeatother

\title{A Twisted Complex Brunn-Minkowski Theorem}
\author{E. M. Ainasse}
\date{\today}
\AtEndDocument{\bigskip{\footnotesize%
  \begin{multicols}{2}
  \textsc{Mathematics Department\\
  Stony Brook University\\
  Stony Brook, NY, 11794-3651, USA}\\
 \vfill\null
 \columnbreak
 \textsc{Department of Applied Mathematics \& Statistics\\
 Stony Brook University\\
 Stony Brook, NY 11794-3600, USA}
  \end{multicols}
   \par  
  \begin{center}
      \textit{E-mail address}: \texttt{elmehdi.ainasse@stonybrook.edu}
  \end{center}}}
  
\begin{document}
\maketitle

\begin{abstract}
In his \textit{Annals of Mathematics} paper \citep{berndtsson-2009-annals}, Berndtsson proves an important result on the Nakano positivity of holomorphic infinite-rank vector bundles whose fibers are Hilbert spaces consisting of holomorphic $L^2$-functions with respect to a family of weight functions $\left\{e^{-\varphi(t,\cdot)}\right\}_{t \in U}$, varying in $t \in U \subset \C^m$, over a pseudoconvex domain. Using a variant of Hörmander’s theorem due
to Donnelly and Fefferman, we show that Berndtsson’s Nakano positivity result holds under different
(in fact, more general) curvature assumptions. This is of particular interest when the manifold admits
a negative non-constant plurisubharmonic function, as these curvature assumptions then allow for some curvature
negativity. We describe this setting as a ``twisted'' setting.
\end{abstract}

\section{Introduction and description of the main results}\label{section:introduction}
Let $X$ be an $n$-dimensional relatively compact complete Kähler submanifold of a Stein Kähler manifold $(Y,g)$. Let $U$ be a domain in $\C^m$ containing the origin, and let $V \rightarrow \overbar{X}$ be a holomorphic vector bundle over the closure $\overbar{X}$ of $X$. Let $\left\{h^{[t]}\right\}_{t \in U}$ be a family of smooth Hermitian metrics for $V \rightarrow \overbar{X}$. By the latter, we mean taking a metric $h$ for the pullback bundle $\pi^{*}_{\overbar{X}} V \rightarrow U \times \overbar{X}$, where $\pi_{\overbar{X}} : U \times \overbar{X} \rightarrow \overbar{X}$ is the canonical projection, and letting $h^{[t]} := i^{*}_t h$ where $i_t$ denotes the inclusion map $\overbar{X} \hookrightarrow \{t\} \times \overbar{X}$. We can consider the space $L^2\left(X,h^{[t]}\right)$ of sections $f$ of $V \rightarrow X$ whose norm $\abs{f}_{h^{[t]}}$, with respect to the metric $h^{[t]}$, is square-integrable on $X$ with respect to the volume form $dV_g$ induced by the Kähler metric $g$. For each $t \in U$, we can then consider the space $\HH\left(X,h^{[t]}\right)$ of holomorphic sections of $V \rightarrow X$ in $L^2\left(X,h^{[t]}\right)$, which is a closed subspace of $L^2\left(X,h^{[t]}\right)$ by our smoothness and boundedness assumptions, and therefore a Hilbert space. The smoothness and boundedness assumptions imply further that the Hilbert spaces in the collection $\left\{\HH\left(X,h^{[t]}\right)\right\}_{t \in U}$ have equivalent norms. Indeed, for any section $f$ of $V \rightarrow X$, for any point $x \in X$, and for any $s,t \in U$,
$$\abs{f(x)}^2_{h^{[t]}} = \sup_{\sigma \in V^{*}_x-\{0\}} \dfrac{\abs{\left\langle\sigma,f(x)\right\rangle}^2}{\abs{\sigma}^2_{h^{[t],*}}} =  \left(\sup_{\sigma \in V^{*}_x-\{0\}} \dfrac{\abs{\sigma}^2_{h^{[s],*}}}{\abs{\sigma}^2_{h^{[t],*}}}\right)\abs{f(x)}^2_{h^{[s]}},$$
where $h^{[t],*}$ and $h^{[s],*}$ are the dual metrics for the dual bundle $V^{*} \rightarrow X$ induced by $h^{[t]}$ and $h^{[s]}$ respectively. Therefore, for any section $f$ of $V \rightarrow X$,
$$C^{-1}_{t,s}\int_X \abs{f(x)}^2_{h^{[s]}}dV_g(x) \leq \int_X \abs{f(x)}^2_{h^{[t]}}dV_g(x) \leq C_{s,t}\int_X \abs{f(x)}^2_{h^{[s]}}dV_g(x),$$
where
$$C_{s,t} := \sup_{\substack{\sigma \in V^{*}_x-\{0\} \\ x \in X}} \dfrac{\abs{\sigma}^2_{h^{[s],*}}}{\abs{\sigma}^2_{h^{[t],*}}},$$
and $C_{t,s}$ is defined similarly. In particular, it follows that for any $s,t \in U$, $f \in \HH\left(X,h^{[t]}\right)$ if and only if $f \in \HH\left(X,h^{[s]}\right)$. Thus, the underlying vector spaces of the Bergman spaces $\mathcal{H}^2\left(X,h^{[t]}\right)$ are equal as subspaces of the space $\holosections(X,V)$. By fixing $\mathcal{H}^2_0 := \HH\left(X,h^{[0]}\right)$, we can define the bundle $E_h$ of infinite rank over $U$ with total space $U \times \mathcal{H}^2_0$, whose fiber over $t \in U$ is $\{t\} \times \mathcal{H}^2_0 \cong \HH_t =: \mathcal{H}^2\left(X,h^{[t]}\right)$. It is a trivial Hilbert bundle equipped with the non-trivial Hermitian metric $\left(\cdot,\cdot\right)_{h^{[t]}}$, varying in $t$, induced by the $L^2$-norm on $\mathcal{H}^2_t$.\\
 
Before stating our main theorem, we need to define what we call a \textit{twisted} curvature operator. Given a smooth Hermitian metric $h$ for $V \rightarrow \overbar{X}$, let $\Theta_{\delta}(h)$ be locally defined as
\begin{align*}
&\sum_{1 \leq j,k \leq m}\dfrac{\del}{\del \bar{t}_k}\left(h^{-1}\dfrac{\del h}{\del t_j}\right) d\bar{t}_k \wedge dt_j + \sum_{\substack{1 \leq j \leq m \\ 1 \leq \mu \leq n}}\dfrac{\del}{\del\bar{z}_\mu}\left(h^{-1}\dfrac{\del h}{\del t_j}\right) d\bar{z}_\mu \wedge dt_j\\
&\, \, \, \, \, \, + \sum_{\substack{1 \leq k \leq m \\ 1 \leq \nu \leq n}}\dfrac{\del}{\del\bar{t}_k}\left(h^{-1}\dfrac{\del h}{\del z_{\nu}}\right) d\bar{t}_k \wedge dz_\nu + \dfrac{\delta}{1+\delta}\sum_{1 \leq \nu, \mu \leq n}\dfrac{\del}{\del\bar{z}_\mu}\left(h^{-1}\dfrac{\del h}{\del z_{\nu}}\right)d\bar{z}_{\mu} \wedge dz_{\nu},
\end{align*}
where $\delta > 0$. We can also express $\Theta_{\delta}(h)$ as
$$\Theta_{\delta}(h) = \Theta(h)-\dfrac{1}{1+\delta}\pi^{*}_X\Theta\left(h^{[t]}\right).$$
As a block matrix split with respect to the product structure $U \times X$ for our holomorphic trivial family, $\Theta_{\delta}(h)$ has the form
$$\Theta_{\delta}(h) = \begin{pmatrix} \dbar_U\left(h^{-1}\del_U h\right) & \dbar_X\left(h^{-1}\del_U h\right) \\ \dbar_U\left(h^{-1}\del_X h\right) & \dfrac{\delta}{1+\delta}\dbar_X\left(h^{-1}\del_X h\right)\end{pmatrix}.$$

Let $\eta$ be a smooth function on $Y$ and define the \textit{twisted} curvature operator $\Xi_{\delta,\eta}$ by
$$\Xi_{\delta,\eta}(h) := \Theta_{\delta}(h) + \dfrac{\delta}{1+\delta}\pi^{*}_X\left(\left(\ricci(g) + 2\ddbar_X\eta-(1+\delta)\del_X\eta\wedge\dbar_X\eta\right)\otimes \mathrm{Id}_V\right).$$

The operator $\Xi_{\delta,\eta}$ can also be represented as the block matrix
$$\Xi_{\delta,\eta}(h) = \begin{pmatrix}
\dbar_U\left(h^{-1}\del_U h\right) & \dbar_X\left(h^{-1}\del_U h\right)\\
\dbar_U\left(h^{-1}\del_X h\right) & \dfrac{\delta}{1+\delta}\left(\dbar_X\left(h^{-1}\del_X h\right) + \ricci(g) + 2\del_X\dbar_X\eta-(1+\delta)\del_X\eta\wedge\dbar_X\eta\right)
\end{pmatrix}.$$
Noting that
$$2\del_X\dbar_X\eta-(1+\delta)\del_X\eta\wedge\dbar_X\eta = \dfrac{4e^{\frac{1+\delta}{2}\eta}}{1+\delta}\del_X\dbar_X\left(-e^{-\frac{1+\delta}{2}\eta}\right),$$ we may also rewrite $\Xi_{\delta,\eta}(h)$ as
$$\Xi_{\delta,\eta}(h) = \begin{pmatrix}
\dbar_U\left(h^{-1}\del_U h\right) & \dbar_X\left(h^{-1}\del_U h\right)\\
\dbar_U\left(h^{-1}\del_X h\right) & \dfrac{\delta}{1+\delta}\left(\dbar_X\left(h^{-1}\del_X h\right) + \ricci(g) + \dfrac{4e^{\frac{1+\delta}{2}\eta}}{1+\delta}\del_X\dbar_X\left(-e^{-\frac{1+\delta}{2}\eta}\right)\right)
\end{pmatrix}.$$

Note that $\Xi_{\delta,\eta}(h) \geq_{\mathrm{Griff}} 0$ if and only if $\Xi_{\delta,\eta}(h) \geq 0$ as a block matrix.\\

Our initial result is stated as follows.

\begin{thmx}\label{thm-nakano-positivity-smooth-bounded}
Let $X$ be an $n$-dimensional relatively compact complete Kähler submanifold of an ambient Stein Kähler manifold $(Y,g)$. Let $V \rightarrow \overbar{X}$ be a holomorphic vector bundle. Let $U \subset \C^m$ be a domain, and let $\left\{h^{[t]}\right\}_{t \in U}$ be a family of smooth Hermitian metrics for $V \rightarrow \overbar{X}$. Let $\delta > 0$ and let $\eta$ be a smooth function on $Y$. If $\Xi_{\delta,\eta}(h) >_{\mathrm{Griff}} 0$ and
$$\dbar_X\left(\left(h^{[t]}\right)^{-1}\del_X h^{[t]}\right) + \left(\ricci(g) + 2\del_X\dbar_X\eta-(1+\delta)\del_X\eta\wedge\dbar_X\eta\right) \otimes \mathrm{Id}_V >_{\mathrm{Nak}} 0,$$
for each $t \in U$, then the holomorphic Hermitian bundle $\left(E_h,\left(\cdot,\cdot\right)_{h^{[t]}}\right)$ is Nakano positive. Moreover, if either $\Xi_{\delta,\eta}(h) \geq_{\mathrm{Griff}} 0$ or $$\dbar_X\left(\left(h^{[t]}\right)^{-1}\del_X h^{[t]}\right) + \left(\ricci(g) + 2\del_X\dbar_X\eta-(1+\delta)\del_X\eta\wedge\dbar_X\eta\right) \otimes \mathrm{Id}_V \geq_{\mathrm{Nak}} 0,$$
then $\left(E_h,\left(\cdot,\cdot\right)_{h^{[t]}}\right)$ is Nakano semipositive.
\end{thmx}

When $V \rightarrow \overbar{X}$ is a line bundle $L \rightarrow \overbar{X}$ equipped with a smooth Hermitian metric $h = e^{-\varphi}$, $\Theta_{\delta}\left(e^{-\varphi}\right)$ can be represented as
$$\Theta_{\delta,\eta}\left(e^{-\varphi}\right) = \begin{pmatrix} \del_U\dbar_U\varphi & \del_U\dbar_X\varphi \\ \del_X\dbar_U\varphi & \dfrac{\delta}{1+\delta}\del_X\dbar_X\varphi\end{pmatrix}$$
and the twisted curvature operator $\Xi_{\delta,\eta}\left(e^{-\varphi}\right)$ can be represented as
$$\Xi_{\delta,\eta}(h) = \begin{pmatrix}
\dbar_U\left(h^{-1}\del_U h\right) & \dbar_X\left(h^{-1}\del_U h\right)\\
\dbar_U\left(h^{-1}\del_X h\right) & \dfrac{\delta}{1+\delta}\left(\dbar_X\left(h^{-1}\del_X h\right) + \ricci(g) + 2\del_X\dbar_X\eta-(1+\delta)\del_X\eta\wedge\dbar_X\eta\right)
\end{pmatrix}.$$

In this case, Griffiths (semi)positivity and Nakano (semi)positivity are equivalent. In particular, the positivity of $\Xi_{\delta,\eta}\left(e^{-\varphi}\right)$ implies the positivity of
\begin{align*}
    &\dbar_X\left(\left(h^{[t]}\right)^{-1}\del_X h^{[t]}\right) + \ricci(g) + 2\del_X\dbar_X\eta-(1+\delta)\del_X\eta\wedge\dbar_X\eta\\
    &= \del_X\dbar_X\varphi^{[t]} + \ricci(g) + 2\del_X\dbar_X\eta-(1+\delta)\del_X\eta\wedge\dbar_X\eta,
\end{align*}
for each $t \in U$, by Schur complement theory, in which case the second curvature condition becomes redundant.\\

Our theorem mainly follows from an application of the Donnelly-Fefferman-Ohsawa theorem on $L^2$-estimates for the $\dbar$-operator. (See \S \ref{section:background-donnelly}.)\\

Theorem \ref{thm-nakano-positivity-smooth-bounded} can be seen as a ``twisted'' variant of Berndtsson's celebrated Nakano positivity theorem \citep[Theorem 1.1]{berndtsson-2009-annals}. If $X = \Omega$ is a bounded pseudoconvex domain in $Y = \C^n$, $V \rightarrow \overbar{X}$ is a trivial bundle $L \rightarrow \overbar{X}$ of rank $1$, and $h = e^{-\varphi}$ is a metric for the pullback bundle $\pi^{*}_{\overbar{X}} L \rightarrow U \times \overbar{X}$, then we are in the context of Berndtsson's Nakano positivity theorem.

\begin{theorem}\label{berndtsson-nakano-positivity-theorem}\emph{(\citep[Theorem 1.1]{berndtsson-2009-annals})}
If $\Omega$ is pseudoconvex and $\varphi$ is plurisubharmonic (resp. strictly plurisubharmonic) on $U \times \Omega$, then the bundle $\left(E,\left(\cdot,\cdot\right)_{\varphi(t,\cdot)}\right)$ is Nakano semipositive (resp. positive).
\end{theorem}

By choosing $\eta$ to be identically constant and letting $\delta \rightarrow +\infty$ in our Theorem \ref{thm-nakano-positivity-smooth-bounded}, we recover Theorem \ref{berndtsson-nakano-positivity-theorem}.\\

Unlike Berndtsson's case in which the curvature hypothesis would be that of Nakano positivity of the metric $h$ for the pullback bundle in this geometric setting (\citep[Theorem 1.5]{raufi2013log}), our curvature hypothesis is more general, and allows for some amount of curvature negativity along the manifold $X$ in certain cases.\\

For simplicity, consider the case where we have a line bundle $L \rightarrow Y$ equipped with a family of smooth Hermitian metrics. We can see from the matrix representation $\Xi_{\delta,\eta}\left(e^{-\varphi}\right)$, if the function $-e^{-\frac{1+\delta}{2}\eta}$ is plurisubharmonic on $X$, then our metric $e^{-\varphi}$ does not need to be positively curved along the fiber $X$ for the condition $\Xi_{\delta,\eta}\left(e^{-\varphi}\right) \geq 0$ to be satisfied. In principle, we may choose a metric $e^{-\varphi}$ so that the curvature along the fiber $X$ is possibly as negative as
$$\dfrac{4e^{\frac{1+\delta}{2}\eta}}{1+\delta}\del_X\dbar_X\left(-e^{-\frac{1+\delta}{2}\eta}\right).$$
Therefore, provided that the manifold $X$ possesses a negative non-constant plurisubharmonic function, our result improves Berndtsson's result. The existence of negative non-constant plurisubharmonic functions is equivalent to the existence of functions \textit{of self-bounded gradient} (see \S \ref{sbg-background}).

\section{Motivating examples}\label{motivating-examples}
\subsection{The unit ball}
Let $\delta = 1$. On the unit ball $\mathbb{B}_n(1)$ in $\C^n$, the function $z \mapsto -\log\left(1-\abs{z}^2\right)$ is a function of self-bounded gradient (with constant $1$) on $\mathbb{B}_n(1)$. Therefore,
$$2\del_z\dbar_z\eta-(1+\delta)\left(\del_z\eta\wedge\dbar_z\eta\right) = \dfrac{4e^{\frac{1+1}{2}\eta}}{1+1}\del_z\dbar_z\left(-e^{-\frac{1+1}{2}\eta}\right) = \dfrac{2}{1-\abs{z}^2}dz \dot{\wedge} d\bar{z},$$
where
$$dz\dot{\wedge} d\bar{z} := dz_1 \wedge d\bar{z}_1 + \cdots + dz_n \wedge d\bar{z}_n.$$

If we equip $\mathbb{B}_n(1)$ with the Euclidean metric, the twisted curvature condition $\Xi_{\delta,\eta}\left(e^{-\varphi}\right) \geq 0$ reduces to the following.
\begin{equation}\label{curvature-unit-ball}
    \begin{pmatrix}
\del_t\dbar_t\varphi & \del_t\dbar_z\varphi \\
\del_z\dbar_t\varphi & \dfrac{1}{2}\left(\del_z\dbar_z\varphi + \dfrac{2}{1-\abs{z}^2}dz\dot{\wedge}d\bar{z}\right)
\end{pmatrix} \geq 0.
\end{equation}

Therefore, $\del_z\dbar_z\varphi$ can be chosen to be as negative as $-\dfrac{2}{1-\abs{z}^2}dz \dot{\wedge} d\bar{z}$.

\begin{example}(A diagonal weight on $\mathbb{B}_n(1)$)
If we pick a weight $\varphi$ of the form $$\varphi(t,z) = \vartheta(t) + \psi(z),$$ then the twisted curvature condition reduces to
$$\begin{pmatrix}
\del_t\dbar_t\vartheta & 0 \\
0 & \dfrac{1}{2}\left(\del_z\dbar_z\varphi + \dfrac{2}{1-\abs{z}^2}dz\dot{\wedge}d\bar{z}\right)
\end{pmatrix} \geq 0,$$
which is equivalent to
$$\del_t\dbar_t\vartheta \geq 0 \text{ and } \del_z\dbar_z\psi + \dfrac{2}{1-\abs{z}^2} dz \dot{\wedge} d\bar{z}\geq 0.$$

For instance, let $\varphi(t,z) = \abs{t}^2 -2\abs{z}^2$. Then
$$\del_t\dbar_t\varphi =  dt \dot{\wedge} d\bar{t} \geq 0 \text{ and } \del_z\dbar_z\psi + \dfrac{2}{1-\abs{z}^2} dz \dot{\wedge} d\bar{z} = \dfrac{\abs{z}^2 }{1-\abs{z}^2} dz \dot{\wedge} d\bar{z} \geq 0.$$

So the consequence of Theorem \ref{berndtsson-nakano-positivity-theorem} still holds, even though $\varphi$ is not plurisubharmonic.
\end{example}

\begin{example} (A non-diagonal weight on $\mathbb{D}$)
Consider now the case when $n = 1$; the unit disk. Suppose that $U$ is a disk with radius $\sqrt{2}$ centered at the origin. Once again, let $\delta = 1$. Another weight one might consider is $$\varphi(t,z) = \left(1-\abs{z}^2\right)\abs{t}^2 = \abs{t}^2-\abs{t}^2\abs{z}^2.$$

Clearly, $\varphi$ is not plurisubharmonic. However, it satisfies condition \eqref{curvature-unit-ball}. Indeed, the trace of the form
\begin{align*}
&\begin{pmatrix}
\I\del_t\dbar_t\varphi & \I\del_t\dbar_z\varphi \\
\I\del_z\dbar_t\varphi & \dfrac{1}{2}\left(\I \del_z\dbar_z\varphi + \dfrac{2}{1-\abs{z}^2}dz \wedge \I d\bar{z}\right)
\end{pmatrix}\\
&= \begin{pmatrix}
(1-\abs{z}^2)  \I dt \wedge d\bar{t} &  -t\bar{z} \I dt \wedge d\bar{z} \\
 -\bar{z}t  \I dz \wedge d\bar{t} & \dfrac{1}{2}\left(-\abs{t}^2+\dfrac{2}{1-\abs{z}^2}\right)  \I dz \wedge d\bar{z}
\end{pmatrix},
\end{align*}
is clearly a positive form. Moreover, its determinant
\begin{align*}
&\left[\dfrac{\delta}{1+\delta}\left(-\abs{t}^2+\dfrac{2}{1-\abs{z}^2}\right)(1-\abs{z}^2)+\abs{t}^2\abs{z}^2\right]\I dt\wedge d\bar{t} \wedge \I dz \wedge d\bar{z}\\
&= \left[\dfrac{\delta}{1+\delta}\left(\abs{t}^2\abs{z}^2+\left(2-\abs{t}^2\right)\right)+\abs{t}^2\abs{z}^2\right]\I dt\wedge d\bar{t} \wedge \I dz \wedge d\bar{z} > 0
\end{align*}
is also a positive form. Thus $\I\Xi_{1,\eta}\left(e^{-\varphi}\right) > 0$.
\end{example}

\subsection{Pseudoconvex domains with smooth boundary in \texorpdfstring{$\C^n$}{}}
Let $\Omega$ be a bounded pseudoconvex domain in $\C^n$ with smooth boundary. Let $K_{\Omega}$ denote the Bergman kernel of $\Omega$. By Fefferman's theorem on the asymptotic expansion of $K_{\Omega}$, (see \citep[Theorem 2]{fefferman1974-2} and \citep{Feffermann1974} for details) the function
$$z \mapsto \eta(z) := \dfrac{1}{n+1}\log\left(K_\Omega(z,\bar{z})\right)$$
satisfies
$$\del_z\dbar_z\eta - \del_z\eta \wedge \dbar_z\eta \geq -C\I\del_z\dbar_z\abs{\cdot}^2$$
for some constant $C > 0$. The result fails if the constant $\dfrac{1}{n+1}$ in $\eta$ is replaced by a larger constant. (See \citep[{\scshape{Theorem 3.7.6}}]{mcneal-varolin} and \citep[pp. 102-103]{varolin2019notes} for details.)\\

If we let $\varphi(t,z) = \abs{t}^2 -(1-\delta)\left(\dfrac{1}{n+1}\log\left(K_\Omega(z,\bar{z})\right)\right)+C(1+\delta)\abs{z}^2$ for $(t,z) \in \C^m \times \Omega$ and $\delta \in (0,1]$, then $\varphi$ fails to be plurisubharmonic near the boundary of $\Omega$ but clearly satisfies the condition $\Xi_{\delta,\eta}\left(e^{-\varphi}\right) \geq 0$. Indeed, near any boundary point, $$\I\del_z\dbar_z\left(\dfrac{1}{n+1}\log\left(K_\Omega(z,\bar{z})\right)\right) = \omega^B_\Omega$$ is arbitrarily close to $\omega^B_{\mathbb{B}_n(r_0)}$, the Bergman metric of a ball of radius $r_0$ centered at the origin in $\C^n$, in suitably chosen coordinates, and the latter is arbitrarily large near any such point in those coordinates.

\section{Background and notation}\label{background}

\subsection{Functions of self-bounded gradient}\label{sbg-background}
\begin{definition}
Let $X$ be a complex manifold. A function $f \in W^{2,1}_{\mathrm{loc}}(X)$ has \textit{self-bounded gradient with constant $c > 0$} is a function such that
$$\ddbar f \geq c\cdot\left(\del f \wedge \dbar f\right).$$
\end{definition}
We denote the set of such functions is denoted by $\SBGc(X)$.\\

\begin{remark}\label{sbgc-sbg1}
For any $c >0$, $f \in \SBGc(X)$ if and only if $cf \in \SBG(X)$, and so we only need to consider functions belonging to $\SBG(X)$. From now on, we refer to any function in $\SBG(X)$ as a \textit{function of self-bounded gradient}. McNeal introduced this term in \citep{mcneal}.\\
\end{remark}

Moreover, $f$ has self-bounded gradient if and only if $-e^{-f}$ is plurisubharmonic, since
$$\ddbar\left(-e^{-f}\right) = e^{-f}\left(\ddbar f -\del f \wedge \dbar f\right).$$

In particular, this means that a complex manifold admits a function of self-bounded gradient if and only if it possesses a negative plurisubharmonic function.\\

Typical examples of functions of self-bounded gradient are the potentials for the Poincaré metric and the hyperbolic metric on the unit ball and the punctured disk, respectively, as seen in \S \ref{section:introduction}. Other interesting examples can be found in \citep{mcneal-varolin}. Additional general examples include the following.

\begin{example}(Strongly pseudoconvex domains) If a domain $\Omega$ is strongly pseudoconvex, then it admits plurisubharmonic defining function $\varrho$. We can then take $\eta := -\log(-\varrho)$ to be our $\SBG$ function.
\end{example}

\begin{example}(Relatively compact pseudoconvex domains with smooth boundary)\label{example-diederich-fornaess}\\
Let $X$ be a Stein manifold and let $\Omega \subset X$ be a relatively compact pseudoconvex subdomain.
\begin{enumerate}
    \item Assuming that the boundary of $\Omega$ is smooth, \citep[{\scshape{Theorem}} 1]{diederich-fornaess} states the existence of a smooth strictly plurisubharmonic function $\rho$ on $\Omega$ with negative values and which converges to zero at the boundary. In this situation, we may choose $\eta := -\log(-\rho)$ to be our function in $\SBG(\Omega)$.
    
    \item If the boundary of $\Omega$ is $\mathcal{C}^{r}$-smooth; $2 \leq r \leq \infty$, another theorem of Diederich-Forn{\ae}ss (\citep[Theorem 1]{diederich-fornaess-2}) states the existence of a defining function $\rho$ that is $\mathcal{C}^r$-smooth in a neighborhood of $\overbar{\Omega}$, and such that $\hat{\rho} := -(-\rho)^{\gamma}$ is a strictly plurisubharmonic bounded exhaustion function on $\Omega$ for any small enough number $\gamma \in (0,1)$. In this situation, we may then choose $\eta = -\log(-\hat{\rho}) = -\gamma\log(-\rho)$ as our function in $\SBG(\Omega)$.
\end{enumerate}
\end{example}

\begin{remark}\label{remark-hed}
The results of Diederich and Forn{\ae}ss have first been extended to relatively compact pseudoconvex domains with $\mathcal{C}^1$ boundary by Kerzman-Rosay \citep{kerzman-rosay}, and then further extended to pseudoconvex domains with Lipschitz boundary by Demailly \citep{Demailly1987}. More recently, Avelin, Hed and Persson extended these results to pseudoconvex domains with Log-Lipschitz boundary (see \citep[Theorem 3, Corollary 4]{hed-benny-persson}).
\end{remark}

\begin{example}(Hyperconvex manifolds)
A \textit{hyperconvex} manifold $X$ is a manifold that admits a bounded strictly plurisubharmonic exhaustion function $\psi : X \rightarrow [-\infty,b)$ (see \citep{stehle}).  Define $\eta := -\log(b-\psi)$. Since $-e^{-\eta} = \psi-b$, then $\eta \in \SBG(X)$.
\end{example}

Clearly, we may assume that the exhaustion is negative.\\

\begin{remark}
Note that not every pseudoconvex domain is hyperconvex. A counterexample is the Hartogs triangle $T = \left\{(z,w) \in \C^2 : \abs{z} < \abs{w} < 1\right\}$ as observed by Diederich and Forn\ae ss in \citep{diederich-fornaess} and \citep{diederich-fornaess-2}.
\end{remark}

\begin{remark}\label{hyperconvex-exhaustion-regularize}
Let $X$ be a complete hyperconvex manifold and let $\psi : X \rightarrow [-\infty,0)$ be its plurisubharmonic exhaustion function. By regularizing $\psi$, $X$ is hyperconvex if and only if it admits a smooth plurisubharmonic exhaustion $\tilde{\psi} : \Omega \rightarrow [-1,0)$.
\end{remark}

\begin{remark}
By Example \ref{example-diederich-fornaess} and Remark \ref{remark-hed}, any relatively compact pseudoconvex with Log-Lipschitz boundary in a Stein manifold is hyperconvex.\\
\end{remark}

Hyperconvex domains and functions of self-bounded gradient are related as follows.

\begin{proposition}\emph{\cite[Proposition 4.6]{Ohsawa2018}}\label{hyperconvex-sbg}
A manifold $X$ is hyperconvex if and only if there exists a strictly plurisubharmonic
exhaustion function $\psi$ on $X$ which is in $\SBGc(X)$ for some positive constant $c$.
\end{proposition}

Although hyperconvex domains can be of interest in other areas of several complex variables (such as pluripotential theory), we will not be directly concerned with them in this article.

\subsection{The Donnelly-Fefferman-Ohsawa Theorem}\label{section:background-donnelly}
Here we state (a version of) the Donnelly-Fefferman-Ohsawa theorem on $L^2$-estimates for the $\dbar$-equation.

\begin{theorem}\label{donnelly-fefferman-ohsawa-l2-estimate-theorem}
Let $X$ be a complete Kähler manifold equipped with a Kähler metric $g$ that is not necessarily complete, and let $V \rightarrow X$ be a holomorphic vector bundle with Hermitian metric $h_0$. Assume there exist a smooth function $\eta$, a positive number $\delta$ and a non-negative $(1,1)$-form $\Phi$ such that
$$\Theta(h_0) + \left(\ricci(g)+2\ddbar\eta-(1+\delta)\del\eta\wedge\dbar\eta\right)\otimes\mathrm{Id}_V \geq_{\mathrm{Nak}} \Phi \otimes \mathrm{Id}_V.$$
Then for every $V$-valued $(0,1)$-form $\alpha$ such that
$$\dbar\alpha = 0 \text{ and } \int_{X}\abs{\alpha}^2_{h_0,\Phi} dV_g < +\infty$$
there exists a measurable section $u$ of $V \rightarrow X$ such that
$$\dbar u = \alpha \text{ and } \int_X \abs{u}^2_{h_0} dV_g \leq \left(\dfrac{1+\delta}{\delta}\right)\int_X\abs{\alpha}^2_{h_0,\Phi}dV_g.$$
\end{theorem}

We refer the reader to \citep{varolin2019notes} for a proof, as well as \citep{mcneal-varolin} for a detailed exposition in the case of line bundles.

\subsection{Holomorphic Hilbert bundles of families of relatively compact complete Kähler submanifolds of Stein manifolds}

\subsubsection{Definitions}\label{trivial-families-smooth-bounded-definitions-section}
Let $\pi_{\overbar{X}}$ denote the projection $U \times \overbar{X} \rightarrow \overbar{X}$ where $U$ is a domain in $\C^m$. Assume that $0 \in U$ without loss of generality. We define a \textit{family $\left\{h^{[t]}\right\}_{t \in U}$ of smooth Hermitian metrics} for $V \rightarrow \overbar{X}$ to be a smooth Hermitian metric $h$ for the pullback bundle $\pi^{*}_{\overbar{X}} V \rightarrow U \times \overbar{X}$. It follows that for each $t \in U$, $h^{[t]} := i_t^{*}h$ is a smooth Hermitian metric, where
$$i_t : V \rightarrow \restr{\pi^{*}_{\overbar{X}} V}{\{t\} \times \overbar{X}}$$
is the natural isomorphism of vector bundles induced by the inclusion of $\overbar{X}$ into the fiber $\{t\} \times \overbar{X}$ of $\pi^{*}_{\overbar{X}} V \rightarrow U \times \overbar{X}$.

In this setting, we can define for each $t \in U$ a Hilbert space
\begin{align*}
    \HH\left(X,h^{[t]}\right) 
    &:= \left\{ f \in \holosections(X,V) : \norm{f}^2_{h^{[t]}} := \int_X \abs{f}^2_{h^{[t]}} dV_g :=  \int_X  h^{[t]}\left(f,f\right) dV_g < +\infty \right\}\\
    & \, = \holosections(X,V) \cap L^2\left(X,h^{[t]}\right),
\end{align*}
where:
$$L^2\left(X,h^{[t]}\right) := \left\{ f \in \Gamma(X,V) : \norm{f}^2_{h^{[t]}} := \int_X \abs{f}^2_{h^{[t]}}dV_g :=  \int_X  h^{[t]}\left(f,f\right) dV_g < +\infty \right\}.$$

Here, $\Gamma(X,V)$ denotes the space of measurable sections of $V \rightarrow X$, $\holosections(X,V)$ denotes the space of holomorphic sections of $V \rightarrow X$, and $dV_g$ denotes the volume form induced by the metric $g$. The norm on $L^2_t$ and its corresponding inner product will be denoted by $\norm{.}_{h^{[t]}}$ and $(\cdot,\cdot)_{h^{[t]}}$ respectively. We can then define the holomorphic Hilbert bundle $E_h$ as the infinite-rank vector bundle -- or Hilbert bundle --
$$U \times \HH_0 \rightarrow U,$$
and we define a Hermitian metric on it by endowing the Hilbert space fiber 
$$\{t\} \times \HH_0 \cong \HH_t$$ 
with the norm $\norm{\cdot}_{h^{[t]}}$.\\

We define the space of sections of $E_h$ as
$$\Gamma\left(E_h\right) := \left\{f\in\Gamma\left(U \times X, \pi^{*}_X V\right) : i^{*}_t\hat{f} \in \HH_t, \forall t \in U\right\}.$$

In particular, the space of holomorphic sections $E_h$ is defined as
$$\holosections\left(E_h\right) := \left\{\hat{f} \in \Gamma\left(E_h\right) : \hat{f} \in \holosections(U \times X, \pi^{*}_X V)\right\}.$$

For $\hat{f} \in \Gamma(E_h)$, we denote $i^{*}_t\hat{f}$ by $\hat{f}^{[t]}$. So all sections are holomorphic on the fibers, and for a section is holomorphic if it is holomorphic in the base variable as well.

We denote by $E^{*}_h$ the dual bundle to $E_h$. This bundle is also trivial. The space of sections of $E^{*}_h$ is defined as
$$\Gamma\left(E^{*}_h\right) := \left\{\xi : E_h \rightarrow \C ; \xi_t := \restr{\xi}{\HH_t} \in \left(\HH_t\right)^{*}\right\}.$$

The bundle $E^{*}_h$ is equipped fiberwise with the non-trivial Hermitian dual norm
$$\norm{\xi}_{*,h^{[t]}} = \sup_{f\in\HH_t-\{0\}}\frac{\abs{\left\langle\xi_t,f\right\rangle}}{\norm{f}_{h^{[t]}}},$$
for each $t \in U$, where $\xi$ is a section of $E^{*}_h$ and $\xi_t := \restr{\xi}{\HH_t}$.

A section $\xi$ is smooth (resp. holomorphic) if for each $\hat{f} \in \Gamma(E_h)$ that is smooth (resp. holomorphic), the function 
$$U \ni t \mapsto \left\langle \xi_t,i^{*}_t \hat{f} \right\rangle \in \C$$
is a smooth (resp. holomorphic) function of $U$.

\subsubsection{Preliminaries}\label{trivial-families-smooth-bounded-preliminaries-section}

Let $F_h$ be the Hilbert bundle whose fiber over $t \in U$ is $L^2_t$ and let $E_h$ be its subbundle whose fiber over $t \in U$ is $\HH_{t}$. Let $\mathcal{P}_t : L^2_t \rightarrow \HH_t$ denote the fiberwise Bergman projection, and let $\mathcal{P}^{\perp}_t$ the fiberwise orthogonal projection of $L^2_t$ onto the orthogonal complement of $\HH_t$. Additionally, denote by $\nabla^{F_h}$ and $\nabla^{E_h}$ the Chern connections of each of $F_h$ and $E_h$; and let $\Theta^{F_h}$ and $\Theta^{E_h}$ denote their respective curvature forms.\\

Choose local coordinates $(z_1, \cdots, z_n)$ for an arbitrary point $z \in X$, and denote by $(t_1, \cdots, t_m)$ the global coordinates $t \in U \subset \C^m$. Let $e_1, \dots, e_r$ be a holomorphic frame for $V \rightarrow X$ and let $H$ denote the local matrix representation of $h$ in this frame, i.e.
$$H = \left(h(e_j,e_k)\right)_{j,k = 1}^r = \begin{pmatrix} h(e_1,e_1) & h(e_1,e_2) & \cdots & h(e_1,e_{r-1}) & h(e_1,e_r) \\
h(e_2,e_1) & \ddots & \ddots & \vdots & \vdots \\
\vdots & \ddots & \ddots & \ddots & \vdots\\
\vdots & \vdots & \ddots & \ddots & h(e_{r-1},e_r)\\
h(e_r,e_1) & h(e_r,e_2) & \cdots & h(e_{r-1},e_r) & h(e_r,e_r)\end{pmatrix}.$$

In addition, for any operator $\mathfrak{d} \in \{\del,\dbar,d\}$ and for any variable $w_i \in \{t_j,\bar{t}_k,z_{\mu},\bar{z}_{\nu}\}$ with $1 \leq j,k \leq m$ and $1 \leq \nu, \mu \leq n$, we let $\mathfrak{d}_{w_i}H$ denote the following matrix in the same holomorphic frame,
$$\mathfrak{d}_{w_i}H := \left(\mathfrak{d}_{w_i}h(e_j,e_k)\right)_{j,k = 1}^r = \begin{pmatrix} \mathfrak{d}_{w_i}h(e_1,e_1) & \mathfrak{d}_{w_i}h(e_1,e_2) & \cdots & \mathfrak{d}_{w_i}h(e_1,e_{r-1}) & \mathfrak{d}_{w_i}h(e_1,e_r) \\
\mathfrak{d}_{w_i}h(e_2,e_1) & \ddots & \ddots & \vdots & \vdots \\
\vdots & \ddots & \ddots & \ddots & \vdots\\
\vdots & \vdots & \ddots & \ddots & \mathfrak{d}_{w_i}h(e_{r-1},e_r)\\
\mathfrak{d}_{w_i}h(e_r,e_1) & \mathfrak{d}_{w_i}h(e_r,e_2) & \cdots & \mathfrak{d}_{w_i}h(e_{r-1},e_r) & \mathfrak{d}_{w_i}h(e_r,e_r)\end{pmatrix}.$$
We adopt the same notation for $h^{[t]}$.\\

Finally, given a section $s$ expressed as $s = \sum_{i = 1}^r s_i e_i$ in this holomorphic frame, for a collection of holomorphic functions $s_1, \dots, s_r$ such that $(s_1, \cdots, s_r) \neq (0, \cdots, 0)$, we represent $s$ by the column vector $S = [s_1 \cdots s_r]^{T}$.\\

From the definition of the Chern connection,
$$d_{t_j}\left(u^{[t]},v^{[t]}\right)_{h^{[t]}} = \left(\nabla^{F_h}_{t_j}u^{[t]},v^{[t]}\right)_{h^{[t]}}+\left(u^{[t]},\nabla^{F_h}_{t_j}v^{[t]}\right)_{h^{[t]}},$$
for any two smooth sections $u$ and $v$ of $F_h$. Therefore, letting $^\dagger$ denote the complex conjugate transpose, we have the following for any two sections $u$ and $v$ of $F_h$.

\begin{align*}
    d_{t_j}\left(u^{[t]},v^{[t]}\right)_{h^{[t]}} &= d_{t_j}\int_{X} h^{[t]}\left(u^{[t]},v^{[t]}\right) dV_g\\
    &= \int_X \del_{t_j}\left[h^{[t]}\left(u^{[t]},v^{[t]}\right) dV_g\right] + \int_X \dbar_{t_j} \left[h^{[t]}\left(u^{[t]},v^{[t]}\right) dV_g\right]\\
    &= \int_X \del_{t_j}\left(\left(V^{[t]}\right)^{\dagger}\right) dV_g + \int_X \dbar_{t_j}\left(\left(V^{[t]}\right)^{\dagger}\right) dV_g\\
    &= \int_X \left(\dbar_{t_j}V^{[t]}\right)^{\dagger}H^{[t]}U^{[t]} dV_g + \int_X \left(V^{[t]}\right)^{\dagger}\dbar_{t_j}H^{[t]}U^{[t]} dV_g + \int_{X} \left(V^{[t]}\right)^{\dagger}H^{[t]}\dbar_{t_j}U^{[t]}dV_g\\
    &\, \, \, \, \, + \int_X \left(\del_{t_j}V^{[t]}\right)^{\dagger}H^{[t]}U^{[t]} dV_g + \int_X V^{[t],\dagger}\del_{t_j}H^{[t]}U^{[t]} dV_g + \int_{X} \left(V^{[t]}\right)^{\dagger}H^{[t]}\del_{t_j}U^{[t]} dV_g.
\end{align*}

Rearranging the terms,
\begin{align*}
d_{t_j}\left(u^{[t]},v^{[t]}\right)_{h^{[t]}} &= \int_X \left[\left(\del_{t_j}+\dbar_{t_j}+\left(H^{[t]}\right)^{-1}\del_{t_j}H^{[t]}\right)V^{[t]}\right]^{\dagger} H^{[t]} U^{[t]} dV_g\\
    &\, \, \, \, \, + \int_X \left(V^{[t]}\right)^{\dagger} H^{[t]} \left[\left(\del_{t_j}+\dbar_{t_j}+\left(H^{[t]}\right)^{-1}\del_{t_j}H^{[t]}\right)U^{[t]}\right] dV_g\\
    &= \int_{X} h^{[t]}\left(\left[d_{t_j}+\left(h^{[t]}\right)^{-1}\del_{t_j}h^{[t]}\right]u^{[t]},v^{[t]}\right) dV_g\\
    &\, \, \, \, \, + \int_{X} h^{[t]}\left(u^{[t]},\left[d_{t_j}+\left(h^{[t]}\right)^{-1}\del_{t_j}h^{[t]}\right]v^{[t]}\right) dV_g\\
    &= \left(\left[d_{t_j}+\left(h^{[t]}\right)^{-1}\del_{t_j}h^{[t]}\right]u^{[t]},v^{[t]}\right)_{h^{[t]}} + \left(u^{[t]},\left[d_{t_j}+\left(h^{[t]}\right)^{-1}\del_{t_j}h^{[t]}\right]u^{[t]}\right)_{h^{[t]}}.
\end{align*}

Clearly, $u^{[t]} \mapsto d_{t_j}u^{[t]}-\left[\left(h^{[t]}\right)^{-1}\del_{t_j}h^{[t]}\right]u^{[t]}$ defines a connection. Moreover, if $u^{[t]}$ is holomorphic in $t$, then $d_{t_j}u^{[t]}-\left[\left(h^{[t]}\right)^{-1}\del_{t_j}h^{[t]}\right]u^{[t]}$ is of bidegree $(1,0)$ and so this operator defines a holomorphic connection. Our previous computation shows that the connection $\nabla^{F_h}$ is metric-compatible, and so it must be the Chern connection of $F_h$. In particular, its $(1,0)$-part is given by $\nabla^{F_h,(1,0)}_{t_j} = \del_{t_j}-\left(h^{[t]}\right)^{-1}\del_{t_j}h^{[t]}$.\\

Thus, the connection form of the Chern connection is given by (wedging with) $\left(h^{[t]}\right)^{-1}\del_{t_j}h^{[t]}$, and so the coefficients of the curvature of $F_h$ are given by $\Theta^{F_h}_{t_j \bar{t}_k} = \dbar_{t_k}\left[\left(h^{[t]}\right)^{-1}\del_{t_j}h^{[t]}\right].$

Therefore, by Griffiths' Curvature Formula (see \citep{berndtsson-2009-annals}),
$$\left(\Theta^{F_h}_{t_j \bar{t}_k}u^{[t]},v^{[t]}\right)_{h^{[t]}} = \left(\mathcal{P}^{\perp}_t\left(\nabla^{F_h}_{t_j}u^{[t]}\right),\mathcal{P}^{\perp}_t\left(\nabla^{F_h}_{t_k}v^{[t]}\right)\right)_{h^{[t]}} + \left(\Theta^{E_h}_{t_j \bar{t}_k}u^{[t]},v^{[t]}\right)_{h^{[t]}}$$ 
for any two smooth sections $u$ and $v$ of $E_h$. Whence if we let $u_1,\dots,u_m$ be any $m$ smooth sections of the bundle $E_h$, then
\begin{equation}\label{curvature-proj}
\sum_{1\leq j,k \leq m} \left(\Theta^{F_h}_{t_j \bar{t}_k} u^{[t]}_j,u^{[t]}_k\right)_{h^{[t]}} = \norm{\mathcal{P}^{\perp}_t\left(\sum_{1 \leq j \leq m}\nabla^{F_h}_{t_j}u^{[t]}_j\right)}_{h^{[t]}}^2 + \sum_{1\leq j,k \leq m} \left(\Theta^{E_h}_{t_j \bar{t}_k} u^{[t]}_j,u^{[t]}_k\right)_{h^{[t]}}.
\end{equation}

The Nakano positivity of $E_h$ will be established by estimating $\sum_{1\leq j,k \leq m}\left(\Theta^{E_h}_{t_j \bar{t}_k}u^{[t]}_j,u^{[t]}_k\right)_{h^{[t]}}$ using the curvature formula (\ref{curvature-proj}). Doing so amounts to estimating the following norm.
\begin{align}\label{perp-identity}
    \norm{\mathcal{P}^{\perp}_t\left(\sum_{1\leq j \leq m}\nabla^{F_h}_{t_j}u^{[t]}_j\right)}^2_{h^{[t]}} &= \norm{\mathcal{P}^{\perp}_t\left(\sum_{1\leq j \leq m}\del_{t_j}u^{[t]}_j-\left(h^{[t]}\right)^{-1}\del_{t_j}h^{[t]}u^{[t]}_j\right)}^2_{h^{[t]}} \nonumber \\ 
    &= \norm{\mathcal{P}^{\perp}_t\left(\sum_{1\leq j \leq m}\left(h^{[t]}\right)^{-1}\del_{t_j}h^{[t]}u^{[t]}_j\right)}^2_{h^{[t]}}.
\end{align}

\section{Proof of Theorem \ref{thm-nakano-positivity-smooth-bounded}}

For the time being, suppose that
$$\dbar_X\left(\left(h^{[t]}\right)^{-1}\del_X h^{[t]}\right) + \left(\ricci(g) + 2\del_X\dbar_X\eta-(1+\delta)\del_X\eta\wedge\dbar_X\eta\right) \otimes \mathrm{Id}_V >_{\mathrm{Nak}} 0,$$
for each $t \in U$, and that $\Xi_{\delta,\eta}\left(h\right) >_\mathrm{Griff} 0$, where $\Xi_{\delta,\eta}$ is the operator introduced in Section \ref{section:introduction}. We will first prove Theorem \ref{thm-nakano-positivity-smooth-bounded} under this assumption. Relaxing our twisted curvature assumptions from positivity to semipositivity is a more delicate process that requires a limiting argument. We show how to do so in Section \ref{trivial-families-smooth-bounded-relax-assumption}.

\subsubsection{Proof of Theorem \ref{thm-nakano-positivity-smooth-bounded} under the assumption of strict curvature positivity}\label{trivial-families-smooth-bounded-proof-strict-section}
By assumption, for each $t \in U$
$$\dbar_X\left(\left(h^{[t]}\right)^{-1}\del_X h^{[t]}\right) + \left(\ricci(g) + 2\del_X\dbar_X\eta-(1+\delta)\del_X\eta\wedge\dbar_X\eta\right) \otimes \mathrm{Id}_V >_{\mathrm{Nak}} 0.$$
By letting
$$\Phi = \dbar_X\left(\left(h^{[t]}\right)^{-1}\del_Xh^{[t]}\right) + \ricci(g) + 2\del_X\dbar_X\eta -(1+\delta)\del_X\eta\wedge\dbar_X\eta,$$
see that the hypotheses of Theorem \ref{donnelly-fefferman-ohsawa-l2-estimate-theorem} are trivially satisfied since $\Theta\left(h^{[t]}\right) = \dbar_X\left(\left(h^{[t]}\right)^{-1}\del_Xh^{[t]}\right)$. Now let $(t_1, \dots, t_m)$ denote the global coordinates of $t \in U$ and let $(z_1, \dots, z_n)$ denote local coordinates for $z \in X$. Consider the section 
$$u^{[t]} := \displaystyle\sum_{1\leq j \leq m}\left(\left(h^{[t]}\right)^{-1}\dfrac{\del h^{[t]}}{\del t_j}\right) u^{[t]}_j.$$ 
Then $u^{[t]}$ solves the $\dbar_z$-equation
$$\dbar_X u^{[t]} =: \dbar_z u^{[t]} = \displaystyle\sum_{\substack{1\leq j \leq m\\ 1 \leq \mu \leq n}}\dfrac{\del}{\del\bar{z}_\mu}\left(\left(h^{[t]}\right)^{-1}\dfrac{\del h^{[t]}}{\del t_j}\right) u^{[t]}_j {d\bar{z}}_{\mu} =: \alpha,$$ 
since every single $u^{[t]}_j$ depends holomorphically on $z = (z_1, \cdots, z_n)$. Note also that
$$\int_X \abs{\alpha}^2_{\Phi,h^{[t]}} dV_g < +\infty.$$

Clearly, the $(0,1)$-form $\alpha$ satisfies $\dbar_X\alpha = 0$. Moreover, since each $u^{[t]}_j$ is in $\HH_t \subset L^2_t$, and since the metric $h$ is smooth up the boundary of $X$, it follows that $u^{[t]}$ is in $\in L^2_t$ as well. By Theorem \ref{donnelly-fefferman-ohsawa-l2-estimate-theorem}, and the fact that $u^{[t]}_0 = \mathcal{P}^{\perp}_t\left(u^{[t]}\right)$ is the minimal-norm solution of $\dbar_X v = \alpha$, we have the estimate
\begin{equation}\label{delta-estimate}
\int_X\abs{u^{[t]}_0}^2_{h^{[t]}} dV_g \leq \left(\dfrac{1+\delta}{\delta}\right)\int_X \abs{\alpha}^2_{\Phi,h^{[t]}}dV_g.    
\end{equation}

Set $\Psi := \Xi_{\delta,\eta}(h)$. Let $\Psi_{ab}$ and $\Psi^{cd}$ denote the components of $\Psi$ in the directions $a$ and $b$, and those of the inverse of $\Psi$ in the directions $c$ and $d$ respectively where $a,b,c,d \in \left\{t_j, \bar{t}_k, z_\mu, \bar{z}_\nu\right\}$. By \eqref{curvature-proj}, \eqref{perp-identity} and \eqref{delta-estimate},

\begin{equation}\label{nakano-estimate}
\sum_{1\leq j,k \leq m} \left(\Theta^{E_h}_{t_j \bar{t}_k}u^{[t]}_j,u^{[t]}_k\right)_{h^{[t]}} \geq \int_X \sum_{1 \leq j,k \leq m}  h^{[t]}\left(\Psi_{t_j \bar{t}_k}u^{[t]}_j,u^{[t]}_k\right) dV_g - \norm{u^{[t]}_0}^2_{h^{[t]}}.
\end{equation}

Furthermore, combining (\ref{nakano-estimate}) with the estimate (\ref{delta-estimate}), we can see that:
\begin{equation}\label{nakano-positive}
\sum_{1\leq j,k \leq m}\left(\Theta^{E_h}_{t_j \bar{t}_k}u^{[t]}_j,u^{[t]}_k\right)_{h^{[t]}} \geq 
\int_X \sum_{1 \leq j,k \leq m} {h^{[t]}}\left(\left[\Psi_{t_j \bar{t}_k} - \sum_{1\leq \mu,\nu \leq n} \Psi^{z_{\mu} \bar{z}_{\nu}}\Psi_{t_j \bar{z}_\mu}\overline{\Psi_{t_k \bar{z}_\nu}}\right]u^{[t]}_j,u^{[t]}_k\right) dV_g.
\end{equation}

Now let $M_{\Psi}$ be the matrix whose $(j,k)$-entries are 
$$\Psi_{t_j \bar{t}_k} - \sum_{1\leq \mu,\nu \leq n} \Psi^{z_{\mu} \bar{z}_{\nu}}\Psi_{t_j \bar{z}_\mu}\overline{\Psi_{t_k \bar{z}_\nu}}.$$

By Schur complement theory,
$\Xi_{\delta,\eta}(h) >_{\mathrm{Griff}} 0$ implies $M_{\Psi} >_{\mathrm{Griff}} 0$. Therefore, we conclude that if $\Xi_{\delta,\eta}(h) >_{\mathrm{Griff}} 0$, then
\begin{equation}\label{nakano-positivity}
\exists c_0 > 0 \text{ (resp. } \geq 0\text{)}: \sum_{1\leq j,k \leq m}\left(\Theta^{E_h}_{t_j \bar{t}_k}u^{[t]}_j,u^{[t]}_k\right)_{h^{[t]}} \geq c_0\sum_{j = 1}^m \norm{u^{[t]}_j}^2_{h^{[t]}}.
\end{equation}

\subsubsection{Relaxing the strict curvature positivity requirement}\label{trivial-families-smooth-bounded-relax-assumption}
Let $h_{\varepsilon} := h e^{-\varepsilon\left(P_{X}^{*}\psi+\abs{t}^2\right)}$ where $\varepsilon > 0$ and $\psi$ is a smooth strictly plurisubharmonic function for the ambient Stein manifold containing $X$, that is smooth up the boundary of $X$. Since $\psi$ is bounded on $X$, we may assume that $\psi > 0$ after subtracting a constant. Moreover, as the result is local, we may assume that $U$ is bounded. Denote by $\Theta^{E_{h_\varepsilon}}_{t_j \bar{t}_k}$ the coefficients for the curvature of $E_{h_{\varepsilon}}$, the Hilbert bundle whose fiber at $t \in U$ is $\HH_{\varepsilon,t} := \HH\left(X,h^{[t]}_{\varepsilon}\right)$. For any $\varepsilon > 0$, the underlying vector spaces of $\mathcal{H}^2_{\varepsilon,t}$ and $\mathcal{H}^2_t$ are equal as subspaces of $\holosections(X,V)$, and so we may act on the same tuple of sections $u_1, \cdots, u_m$.

By construction,
$$\Xi_{\delta,\eta}\left(h_{\varepsilon}\right) = \Xi_{\delta,\eta}(h) + \varepsilon\begin{pmatrix}
\abs{t}^2 & 0 \\ 0 & \del_X\dbar_X\psi 
\end{pmatrix},$$
and
\begin{align*}
&\dbar_X\left(\left(h^{[t]}_{\varepsilon}\right)^{-1}\del_X h^{[t]}_{\varepsilon}\right) + \left(\ricci(g) + 2\del_X\dbar_X\eta-(1+\delta)\del_X\eta\wedge\dbar_X\eta\right) \otimes \mathrm{Id}_V\\
&= \dbar_X\left(\left(h^{[t]}\right)^{-1}\del_X h^{[t]}\right) + \left(\ricci(g) + 2\del_X\dbar_X\eta-(1+\delta)\del_X\eta\wedge\dbar_X\eta\right) \otimes \mathrm{Id}_V + \varepsilon\left(\del_X\dbar_X\psi\otimes\mathrm{Id}_V\right).
\end{align*}
Since
$$\begin{pmatrix}
\abs{t}^2 & 0 \\ 0 & \del_X\dbar_X\psi 
\end{pmatrix} >_{\mathrm{Griff}} 0 \text{ and } \del_X\dbar_X\psi\otimes\mathrm{Id}_V >_{\mathrm{Nak}} 0,$$
it follows that
$$\Xi_{\delta,\eta}\left(h_{\varepsilon}\right) >_{\mathrm{Griff}} \Xi_{\delta,\eta}(h)$$
and 
\begin{align*}
    &\dbar_X\left(\left(h^{[t]}_{\varepsilon}\right)^{-1}\del_X h^{[t]}_{\varepsilon}\right) + \left(\ricci(g) + 2\del_X\dbar_X\eta-(1+\delta)\del_X\eta\wedge\dbar_X\eta\right) \otimes \mathrm{Id}_V\\ &>_{\mathrm{Nak}} \dbar_X\left(\left(h^{[t]}\right)^{-1}\del_X h^{[t]}\right) + \left(\ricci(g) + 2\del_X\dbar_X\eta-(1+\delta)\del_X\eta\wedge\dbar_X\eta\right) \otimes \mathrm{Id}_V.
\end{align*}

Therefore, if either $\Xi_{\delta,\eta}(h) \geq_{\mathrm{Griff}} 0$ or
$$\dbar_X\left(\left(h^{[t]}\right)^{-1}\del_X h^{[t]}\right) + \left(\ricci(g) + 2\del_X\dbar_X\eta-(1+\delta)\del_X\eta\wedge\dbar_X\eta\right) \otimes \mathrm{Id}_V \geq_{\mathrm{Nak}} 0,$$

for each $t \in U$, then
\begin{equation}\label{nakano-positivity-epsilon}
\forall \varepsilon > 0, \exists c^{(\varepsilon)}_0 > 0: \sum_{1\leq j,k \leq m}\left(\Theta^{E_{h_\varepsilon}}_{t_j \bar{t}_k}u^{[t]}_j,u^{[t]}_k\right)_{h^{[t]}_{\varepsilon}} \geq c^{(\varepsilon)}_0\sum_{j = 1}^m \norm{u^{[t]}_j}^2_{h^{[t]}} \geq 0.
\end{equation}

Let $\mathcal{P}^{\varepsilon,\perp}_t$ denote the orthogonal projection of $L^2_{\varepsilon,t}$ onto $\left(\HH_{\varepsilon,t}\right)^{\perp}$. As before, let $\Psi := \Xi_{\delta,\eta}(h)$, and let $\Psi_{ab}$ and $\Psi^{cd}$ denote the components of $\Psi$ in the directions $a$ and $b$, and those of the inverse of $\Psi$ in the directions $c$ and $d$ respectively where $a,b,c,d \in \left\{t_j, \bar{t}_k, z_\mu, \bar{z}_\nu\right\}$. Let us also adopt the same notation for $h_{\varepsilon}$ by letting $\Psi^{\varepsilon}$ represent the corresponding matrix for $h_{\varepsilon}$. Furthermore, let $\Upsilon$ and $\Upsilon^{\varepsilon}$ denote the matrices corresponding to $h^{-1}\del h$ and $h^{-1}_\varepsilon\del h_{\varepsilon}$, respectively.\\

By Griffiths' Curvature Formula, and noting that $\psi$ is independent of $t$,
\begin{align*}
\sum_{1\leq j,k \leq m} \left(\Theta^{E_{h_\varepsilon}}_{t_j \bar{t}_k}u^{[t]}_j,u^{[t]}_k\right)_{h^{[t]}_{\varepsilon}} &= \sum_{1\leq j,k \leq m}\left(\Psi^{\varepsilon}_{t_j\bar{t}_k}u^{[t]}_j,u^{[t]}_k\right)_{h^{[t]}_{\varepsilon}}-\norm{\mathcal{P}^{\varepsilon,\perp}_t\left(\sum_{1\leq j \leq m}\Upsilon^{\varepsilon}_{t_j}u^{[t]}_j\right)}^2_{h^{[t]}_{\varepsilon}}\\
&=\sum_{1\leq j,k \leq m}\left[\left(\Psi_{t_j\bar{t}_k}u^{[t]}_j,u^{[t]}_k\right)_{h^{[t]}_{\varepsilon}} + \varepsilon\left((\pi^{*}_{\overbar{X}} \psi + \abs{t}^2)_{t_j \bar{t}_k}u^{[t]}_j,u^{[t]}_k\right)_{h^{[t]}_{\varepsilon}}\right]\\
&\ \ \ \ \ \ - \norm{\mathcal{P}^{\varepsilon,\perp}_t\left(\sum_{1\leq j\leq m}\Upsilon_{t_j}u^{[t]}_j\right) - \varepsilon\cdot\mathcal{P}^{\varepsilon,\perp}_t\left(\sum_{1\leq j\leq m}(\pi^{*}_{\overbar{X}} \psi + \abs{t}^2)_{t_j}u^{[t]}_j\right)}^2_{h^{[t]}_{\varepsilon}}\\
&= \sum_{1\leq j,k \leq m}\left(e^{-\varepsilon\left(\psi+\abs{t}^2\right)}\Psi_{t_j\bar{t}_k}u^{[t]}_j,u^{[t]}_k\right)_{h^{[t]}} + \varepsilon\sum_{1 \leq j,k \leq m} \delta_{j\bar{k}}\left(e^{-\varepsilon\left(\psi+\abs{t}^2\right)}u^{[t]}_j,u^{[t]}_k\right)_{h^{[t]}}\\
&\ \ \ \ \ \ -\norm{e^{-\varepsilon\left(\psi+\abs{t}^2\right)/2}\mathcal{P}^{\varepsilon,\perp}_t\left(\sum_{1\leq j \leq m}\left(\Upsilon_{t_j}-\varepsilon\bar{t}_j\right)u^{[t]}_j\right)}^2_{h^{[t]}}.
\end{align*}

By adding and subtracting $\sum_{1\leq j,k \leq m}\left(\Theta^{E_h}_{t_j \bar{t}_k}u^{[t]}_j,u^{[t]}_k\right)_{h^{[t]}}$, we then have:
\begin{equation}
\sum_{1\leq j,k \leq m} \left(\Theta^{E_{h_\varepsilon}}_{t_j \bar{t}_k}u^{[t]}_j,u^{[t]}_k\right)_{h^{[t]}_{\varepsilon}} = \sum_{1\leq j,k \leq m}\left(\Theta^{E_h}_{t_j \bar{t}_k}u^{[t]}_j,u^{[t]}_k\right)_{h^{[t]}} + \mathfrak{R}(\varepsilon),
\end{equation}
where
\begin{align*}
\mathfrak{R}(\varepsilon) := & \sum_{1\leq j,k \leq m}\left(\left(e^{-\varepsilon\left(\psi+\abs{t}^2\right)}-1\right)\Psi_{t_j\bar{t}_k}u^{[t]}_j,u^{[t]}_k\right)_{h^{[t]}} + \varepsilon\sum_{1 \leq j,k \leq m} \delta_{j\bar{k}}\left(e^{-\varepsilon\left(\psi+\abs{t}^2\right)}u^{[t]}_j,u^{[t]}_k\right)_{h^{[t]}}\\
&\ \ \ +\norm{\mathcal{P}^{\perp}_t\left(\sum_{1\leq j \leq m}\Upsilon_{t_j}u^{[t]}_j\right)}^2_{h^{[t]}}-\norm{e^{-\varepsilon\left(\psi+\abs{t}^2\right)/2}\mathcal{P}^{\varepsilon,\perp}_t\left(\sum_{1\leq j \leq m}\Upsilon_{t_j}u^{[t]}_j\right)}^2_{h^{[t]}}\\
&\ \ \ +2\varepsilon\mathrm{Re}\left[\left(\mathcal{P}^{\varepsilon,\perp}_t\left(\sum_{1 \leq j \leq m}\Upsilon_{t_j} u^{[t]}_j\right),e^{-\left(\psi+\abs{t}^2\right)}\mathcal{P}^{\varepsilon,\perp}_t\left(\sum_{1 \leq j \leq m}\bar{t}_j u^{[t]}_j\right)\right)_{h^{[t]}}\right]\\
&\ \ \ -\varepsilon^2\norm{e^{-\left(\psi+\abs{t}^2\right)/2}\mathcal{P}^{\varepsilon,\perp}_t\left(\sum_{1 \leq j \leq m}\bar{t}_j u^{[t]}_j\right)}^2_{h^{[t]}}
\end{align*}

In particular,
\begin{equation}\label{epsilon-equality}
\sum_{1\leq j,k \leq m}\left(\Theta^{E_h}_{t_j \bar{t}_k}u^{[t]}_j,u^{[t]}_k\right)_{h^{[t]}} \geq -\mathfrak{R}(\varepsilon),
\end{equation}

Since $u^{[t]}_k \in \mathcal{H}^2_t$ for each $k$ and for each $t \in U$, the first two summands in $\mathfrak{R}(\varepsilon)$ converge to $0$ as $\varepsilon \rightarrow 0$ by smoothness, boundedness, countinuity, and the Cauchy-Schwarz inequality. Since $\psi > 0$ by assumption and $\mathcal{P}^{\varepsilon,\perp}_t$ is an orthogonal projection,
\begin{equation}\label{bound-1}
    \norm{e^{-\left(\psi+\abs{t}^2\right)}\mathcal{P}^{\varepsilon,\perp}_t\left(\sum_{1 \leq j \leq m}\bar{t}_j u^{[t]}_j\right)}^2_{h^{[t]}} \leq \norm{\mathcal{P}^{\varepsilon,\perp}_t\left(\sum_{1 \leq j \leq m}\bar{t}_j u^{[t]}_j\right)}^2_{h^{[t]}} \norm{\sum_{1 \leq j \leq m}\bar{t}_j u^{[t]}_j}^2_{h^{[t]}},
\end{equation}
and
\begin{equation}\label{bound-2}
    \norm{\mathcal{P}^{\varepsilon,\perp}_t\left(\sum_{1 \leq j \leq m}\Upsilon_{t_j} u^{[t]}_j\right)}^2_{h^{[t]}} \leq \norm{\sum_{1 \leq j \leq m}\Upsilon_{t_j} u^{[t]}_j}^2_{h^{[t]}}.
\end{equation}

Each of the upper bounds in \eqref{bound-1} and \eqref{bound-2} respectively is finite by smoothness, boundedness, and the fact that $u^{[t]}_k \in \mathcal{H}^2_t$ for each $1 \leq k \leq m$ and for each $t \in U$. Therefore, the last summand in $\mathfrak{R}(\varepsilon)$ converge to $0$ as $\varepsilon \rightarrow 0$, as does the fifth summand by the Cauchy-Schwarz inequality.\\

It now remains to estimate the difference term $\mathfrak{R}(\varepsilon)$. Note that for any section $\mathfrak{u}$ of $E_h$ and any point $w \in X$,
$$\mathcal{P}^{\varepsilon,\perp}_t \mathfrak{u}^{[t]}(w) = \mathcal{P}^{\perp}_t \mathfrak{u}^{[t]}(w) + \left(\mathfrak{u}^{[t]},K_t(\cdot,w)-e^{-\varepsilon\left(\psi(w)+\abs{t}^2\right)}K^{\varepsilon}_t(\cdot,w)\right)_{h^{[t]}},$$
where $K_t$ and $K^{\varepsilon}_t$ denote the Bergman kernels for $\HH_t$ and $\mathcal{H}^2_{\varepsilon,t}$ respectively. (The emphasis on $w$ is to indicate that the second variable in each Bergman kernel is fixed.) Thus
\begin{align*}
&\mathcal{P}^{\perp}_t \mathfrak{u}^{[t]}(w) - e^{-\varepsilon\left(\psi(w)+\abs{t}^2\right)/2}\cdot\mathcal{P}^{\varepsilon,\perp}_t \mathfrak{u}^{[t]}(w)\\
&= (1-e^{-\varepsilon\left(\psi(w)+\abs{t}^2\right)/2})\cdot\mathcal{P}^{\perp}_t \mathfrak{u}^{[t]}(w) + e^{-\varepsilon\left(\psi(w)+\abs{t}^2\right)/2}\cdot\left(\mathfrak{u}^{[t]},e^{-\varepsilon\left(\psi(w)+\abs{t}^2\right)}K^{\varepsilon}_t(\cdot,w)-K_t(\cdot,w)\right)_{h^{[t]}},
\end{align*}
and so we have the following estimate.
\begin{align*}
&\norm{\mathcal{P}^{\perp}_t \mathfrak{u}^{[t]} - e^{-\varepsilon\left(\psi+\abs{t}^2\right)/2}\cdot\mathcal{P}^{\varepsilon,\perp}_t \mathfrak{u}^{[t]}}_{L^{\infty}_t}\\
&\leq M(\varepsilon)\norm{\mathcal{P}^{\perp}_t  \mathfrak{u}^{[t]}}_{L^{\infty}_t} + e^{-\varepsilon m_0/2}\norm{\mathfrak{u}^{[t]}}_{L^{\infty}_t} m(\varepsilon)\norm{K^{\varepsilon}_t(\cdot,w)}_{L^{\infty}_t}\\
&\, \, \, \, \, + e^{-\varepsilon m_0/2}\norm{\mathfrak{u}^{[t]}}_{L^{\infty}_t} \norm{K^{\varepsilon}_t(\cdot,w)-K_t(\cdot,w)}_{L^{\infty}_t},
\end{align*}
where $M(\varepsilon) = \max\left(\abs{1-e^{-\varepsilon (M_0+R_0)/2}},\abs{1-e^{-\varepsilon m_0/2}}\right), m(\varepsilon) = \max\left(\abs{e^{-\varepsilon (M_0+R_0)}-1},\abs{e^{-\varepsilon m_0}-1}\right)$ and $m_0$ and $M_0$ are the minimum and maximum of $\psi$ over $\overbar{X}$, respectively, and $R_0 := \sup_{\overbar{U}}\abs{t}^2$.

Let $\mathfrak{u} := \sum_{1 \leq j \leq m}\Upsilon_{t_j}u_j$. By smoothness and boundedness, 
$$\norm{\mathfrak{u}^{[t]}}_{L^{\infty}_t}, \norm{\mathcal{P}^{\perp}_t  \mathfrak{u}^{[t]}}_{L^{\infty}_t}, \norm{K^{\varepsilon}_t(\cdot,w)}_{L^{\infty}_t},\norm{K_t(\cdot,w)}_{L^{\infty}_t} < +\infty.$$
Moreover, since $\psi > 0$ by assumption, the sequence of metrics $\left\{h^{[t]}_\varepsilon\right\}_{\varepsilon > 0}$ increases to $h^{[t]}$. Therefore, by the generalization of Ramadanov's theorem on Bergman kernels in \citep{pasternak-winiarski-wojcicki-weighted-ramadanov},
$$\norm{K^{\varepsilon}_t(\cdot,w)-K_t(\cdot,w)}_{L^{\infty}_t} \xrightarrow[\varepsilon \rightarrow 0]{} 0.$$
Finally, knowing that $M(\varepsilon), m(\varepsilon) \rightarrow 0$ as $\varepsilon \rightarrow 0$, we see that 
$$\norm{\mathcal{P}^{\perp}_t \left(\sum_{1\leq j \leq m}\Upsilon_{t_j}u^{[t]}_j\right) - e^{-\varepsilon\left(\psi+\abs{t}^2\right)/2}\cdot\mathcal{P}^{\varepsilon,\perp}_t \left(\sum_{1\leq j \leq m}\Upsilon_{t_j}u^{[t]}_j\right)}_{L^{\infty}_t} \xrightarrow[\varepsilon \rightarrow 0]{} 0,$$ 
and so 
$$\norm{\mathcal{P}^{\perp}_t \left(\sum_{1\leq j \leq m}\Upsilon_{t_j}u^{[t]}_j\right) - e^{-\varepsilon\left(\psi+\abs{t}^2\right)/2}\cdot\mathcal{P}^{\varepsilon,\perp}_t \left(\sum_{1\leq j \leq m}\Upsilon_{t_j}u^{[t]}_j\right)}_{L^{2}_t} \xrightarrow[\varepsilon \rightarrow 0]{} 0.$$ 

Thus,
$$\norm{\mathcal{P}^{\perp}_t \left(\sum_{1\leq j \leq m}\Upsilon_{t_j}u^{[t]}_j\right)}^2_{h^{[t]}} - \norm{e^{- \varepsilon\left(\psi+\abs{t}^2\right)/2}\cdot\mathcal{P}^{\varepsilon,\perp}_t \left(\sum_{1\leq j \leq m}\Upsilon_{t_j}u^{[t]}_j\right)}^2_{h^{[t]}} \xrightarrow[\varepsilon \rightarrow 0]{} 0.$$

\section*{Acknowledgements}
The author would like to express his sincere gratitude to his doctoral advisor Prof. Dror Varolin for continued guidance up until the completion of this article. The author would also like to thank Prof. Eric Bedford, Prof. Christian Schnell, Prof. Gideon Maschler, Mohammed El Alami, Pranav Upadrashta and John Sheridan for inspiring discussions.

\printbibliography
\lipsum[187-201]
\end{document}